\documentclass[letterpaper, 10 pt, conference]{ieeeconf}  

\IEEEoverridecommandlockouts                              

\usepackage{graphicx} 
\usepackage{amsmath} 
\usepackage{amssymb}  
\usepackage{float}
\usepackage{booktabs}
\usepackage{accents}
\usepackage[ruled,vlined]{algorithm2e}
\usepackage{subcaption}
\usepackage{mathtools}
\usepackage{textcomp}
{
    
}

\newcommand{\ignore}[1]{}

\newcommand{\bma}[1]{\left[\begin{array}{#1}}
\newcommand{\ema}{\end{array}\right]}

\DeclareMathAlphabet{\mbf}{OT1}{ptm}{b}{n}

\def\fdotb{{\raisebox{-0.6ex}{ \kern0.2ex\raisebox{0.8ex}{\tiny $\hspace*{-1ex}\circ$}}}}
\def\fddotb{{\raisebox{-0.6ex}{ \kern0.2ex\raisebox{0.8ex}{\tiny $\hspace*{-1ex}\circ\circ$}}}}

\newcommand{\utimes}{ {\raisebox{-0.6ex}{ \kern-1.0ex\raisebox{0.6ex}{ \small $\mathsf{v}$}}} } %
\newcommand{\beq}{\begin{equation}}
\newcommand{\eeq}{\end{equation}}
\newcommand{\bdis}{\begin{displaymath}}
\newcommand{\edis}{\end{displaymath}}
\newcommand{\beqarray}{\begin{eqnarray}}
\newcommand{\eeqarray}{\end{eqnarray}}
\newcommand{\beqarraynn}{\begin{eqnarray*}}
\newcommand{\eeqarraynn}{\end{eqnarray*}}

\newtheorem{assum}{Assumption}
\newtheorem{thrm}{Theorem}

\DeclareMathOperator*{\argmax}{arg\,max}
\DeclareMathOperator*{\argmin}{arg\,min}

\SetKw{Break}{break}

\title{\LARGE \bf State-Dependent Uncertainty Modeling in Robust Optimal Control Problems through Generalized Semi-Infinite Programming}

\author{Jad Wehbeh$^{1}$ and Eric C. Kerrigan$^{2}$
\thanks{$^{1}$Jad Wehbeh is with the Department of Electrical and Electronic Engineering, Imperial College London, SW7 2AZ, UK
        {\tt\small j.wehbeh22@imperial.ac.uk}}%
\thanks{$^{2}$Eric C. Kerrigan is with the Department of Electrical and Electronic Engineering and the Department of Aeronautics, Imperial College London,
        SW7 2AZ, UK
        {\tt\small e.kerrigan@imperial.ac.uk}}
}

\begin{document}

\maketitle
\thispagestyle{empty}
\pagestyle{empty}

\begin{abstract}
Generalized semi-infinite programs (generalized SIPs) are problems featuring a finite number of decision variables but an infinite number of constraints. They differ from standard SIPs in that their constraint set itself depends on the choice of the decision variable. Generalized SIPs can be used to model robust optimal control problems where the uncertainty itself is a function of the state or control input, allowing for a less conservative alternative to assuming a uniform uncertainty set over the entire decision space. In this work, we demonstrate how any generalized SIP can be converted to an existence-constrained SIP through a reformulation of the constraints and solved using a local reduction approach, which approximates the infinite constraint set by a finite number of scenarios. This transformation is then exploited to solve nonlinear robust optimal control problems with state-dependent uncertainties. We showcase our proposed approach on a planar quadrotor simulation where it recovers the true generalized SIP solution and outperforms a SIP-based approach with uniform uncertainty bounds.  
\end{abstract}


\section{Introduction}

\subsection{Background and Motivation}

From the operation of self-driving vehicles to the automation of industrial processes, the current problems of interest in control theory generally focus on dealing with a common challenge: uncertainty~\cite{annaswamy2024control}. Whether this arises from imperfections in the modeling of the system, sensor noise, or other external disturbances, uncertainty in the system dynamics can make it difficult to design controllers that guarantee stability or optimize performance~\cite{weinmann2012uncertain}. Numerous methods have been proposed to address this issue, differing in how they model the uncertainty, the performance guarantees they provide, and the complexity of their implementation.

Several approaches, including many learning-based and some adaptive control methods, seek to avoid the explicit modeling of uncertainty by inferring its distribution from data acquired during system operation~\cite{brunke2022safe}. These methods can simplify the implementation, particularly in complex, high-dimensional systems, and may offer good empirical performance. However, they present challenges in providing performance guarantees, which can be problematic in safety-critical applications where predictable behavior is crucial. On the other hand, stochastic control approaches explicitly model uncertainty using probability distributions and aim to optimize the expected performance while maximizing the likelihood of satisfying system constraints ~\cite{aolaritei2023capture}. While these methods perform well when the uncertainty distributions are accurately known, they often lack guarantees on worst-case performance, and quantifying the uncertainty distributions themselves can be a difficult and resource-intensive task.

In this paper, we specifically focus on robust optimal control methods---a specialized class of robust control techniques that not only guarantee constraint satisfaction under bounded uncertainties but also achieve the best possible performance in the worst-case scenario of these uncertainties~\cite{zagorowska2024automatic}. Building on our previous work on robust optimal control for continuous nonlinear systems~\cite{wehbeh2024semi}, we now address cases where the uncertainty bounds are functions of the state or control input, making them dependent on the decision variables in the optimization problem. This extension is significant because it accommodates more complex uncertainty models, broadening the range of problems that can be effectively tackled while providing a less conservative alternative to assuming the maximum uncertainty bounds over the entire decision space. Problems with state- or control-dependent uncertainties can include examples from aviation \cite{chwa2014fuzzy}, robotics~\cite{navvabi2019new}, electronics \cite{zhao2013robust}, and many other applications.

\subsection{State of the Art}

Explicit methods for handling state-dependent uncertainty are relatively uncommon, with most existing approaches relying on learning-based techniques. For example, \cite{soloperto2018learning} proposes a learning-based robust model predictive control (MPC) framework for linear systems, where both the system and uncertainty models are learned from data using Gaussian process regression. Similarly, \cite{bonzanini2021learning} develops a learning-based stochastic MPC approach that allows uncertainties to depend on the state. Meanwhile, \cite{li2022adaptive} presents an adaptive sliding mode control framework that addresses state-dependent uncertainty through adaptive techniques, without the need for learning-based models. However, none of these approaches provide strict guarantees on the robustness of the system, which makes them inadequate for our application.

One relevant approach to solving robust optimal control problems is through a transformation into a semi-infinite program  (SIP)~\cite{zagorowska2024automatic}. The robust constraint is formulated as an infinite constraint over the compact uncertainty space, and the optimal robust control input is found by solving the SIP through any appropriate method~\cite{djelassi2021recent}. The work of \cite{puschke2018robust}, for instance, explains how to robustly satisfy path constraints in the optimization of batch processes by solving an appropriate SIP. Problems with state-dependent uncertainties, however, are more naturally cast as generalized SIPs, or SIPs where the uncertainty set depends on the decision parameter \cite{stein2003solving}. Generalized SIPs are significantly harder to handle in practice than regular SIPs, though several approaches to solving them have been proposed \cite{vazquez2008generalized}. 

We are specifically interested in approaches that transform the generalized SIPs into standard SIPs because they allow us to use the local reduction framework previously developed in~\cite{wehbeh2024semi}. Originally described in \cite{blankenship1976infinitely}, local reduction discretizes the infinite constraint set into a finite number of scenarios that are obtained through iterative solutions of optimization problems that maximize the constraint violation under the current best guess for the decision variables. The method exhibits good performance in practice and scales well to the high-dimensional problems typically seen in predictive control \cite{hettich1993semi}. Unfortunately, converting generalized SIPs into SIPs is not trivial, and existing methods either rely on restrictive assumptions to define the transformations globally or relax these assumptions to build more limited local transformations \cite{still1999generalized}. The simpler transformation introduced in this work converts generalized SIPs to existence-constrained SIPs, which can be tackled using our previous work in \cite{wehbeh2024semi} or using the global optimization method of \cite{djelassi2021global}.

\subsection{Contributions}

In this paper, we demonstrate how robust optimal control problems with state-dependent uncertainties can be described by generalized SIPs, and present a novel approach to solving such problems. In order to do so, we introduce a method for the conversion of generalized SIPs to existence-constrained SIPs that do not feature any dependence between the uncertainty constraint space and the decision variables. We then show how our extension to the local reduction method of~\cite{blankenship1976infinitely} presented in \cite{wehbeh2024semi} can be used to solve any generalized SIP, and illustrate the mechanics of using this framework for the solution of the robust optimal control problems of interest. We also present some simulation results on a planar quadrotor that showcase the performance of our proposed approach. The principle novelties in the paper are:
\begin{itemize}
    \item The conversion of any generalized SIP to an existence-constrained SIP through the smoothing of logical constraints (Theorem \ref{thrm:transformation} in this paper).
    \item The use of this transformation within the context of our previous work \cite{wehbeh2024semi} to model state-dependent uncertainties in nonlinear robust optimal control problems (Theorem \ref{thrm:convergence} in this paper).  
    \item The demonstration of a practical method for the solution of nonlinear robust optimal control problems with state-dependent uncertainties through simulation results.
\end{itemize}


\section{System Description}

Let $x$ $\in$ $\mathcal{X}$ $\subseteq$ $\mathbb{R}^{n_x}$ be the state trajectory across the period of interest of the discrete-time system with dynamics
\begin{equation}
    \label{eq:dynamics}
    f(x,u,w) = 0
\end{equation}
where  $u$ $\in$ $\mathcal{U}$ $\subseteq$ $\mathbb{R}^{n_u}$ is the collection of control inputs, feedback gains, or other decision variables over this time period, $w$ $\in$ $\mathcal{W}(x,u)$ $\subseteq$ $\mathbb{R}^{n_w}$ represents any unknowns or disturbances affecting the state evolution of the dynamics, and $f(\cdot,\cdot,\cdot):$ $\mathcal{X} \times \mathcal{U} \times \mathcal{W} \rightarrow \mathbb{R}^{n_z}$. $f$ includes all of the dynamics across the time steps considered and is a condensed formulation of the dynamics across successive time steps that is more commonly seen. 

Here, the set of possible uncertainties $\mathcal{W}(x,u)$ can depend on both the state and control trajectories, and is defined as 
\begin{equation}
\label{eq:uncertainty_set}
\mathcal{W}(x,u) \coloneqq \left\{ w \left| \:
    h(x,u,w) \geq 0
    \right. \right\}
\end{equation}
where $h(\cdot,\cdot,\cdot):$ $\mathcal{X} \times \mathcal{U} \times \mathcal{W} \rightarrow \mathbb{R}^{n_h}$.

\begin{assum}
    \label{assum:f_uniqueness}
    The dynamics of $f$ are unique such that for any choice of $u \in \mathcal{U}$ and $w \in \mathbb{R}^{n_w}$, there exists a unique solution $x = \bar{\mathcal{X}}(u,w) \in \mathcal{X}$ to \eqref{eq:dynamics}. Additionally, for any choice of $x \in \mathcal{X}$ and $u \in \mathcal{U}$, the solution space for $w$ to the equation $f(x,u,w) = 0$ is compact.
\end{assum}

This allows us to define the superset of \eqref{eq:uncertainty_set},
\begin{equation}
\mathcal{Y}(u) \coloneqq \left\{ w \left| \:
    h(\bar{\mathcal{X}}(u,w),u,w) \geq 0
    \right. \right\}
\end{equation}
which includes all of the uncertainties that are admissible given a control trajectory $u$. Assuming feasible trajectories for $x$, $\mathcal{W}(x,u) \subset \mathcal{Y}(u)$, since $\mathcal{Y}(u)$ covers the uncertainties that are compatible with any of the trajectories for $x$ that can arise from a choice of $u$, while $\mathcal{W}(x,u)$ is limited to a single one.

We then consider the cost function $J(\cdot,\cdot,\cdot):$ $\mathcal{X} \times \mathcal{U} \times \mathcal{W} \rightarrow \mathbb{R}$, which maps choices of $x$, $u$, and $w$ onto some objective we seek to minimize. We also consider its worst-case realization for a given choice of $u$, $J^*(\cdot):$ $\mathcal{U} \rightarrow \mathbb{R}$, which can be described as
\begin{subequations}
\begin{equation}
    J^*(u) \coloneqq \max_{w  \,\in \, \mathcal{Y}(u)} \, J(x,u,w)
\end{equation}
s.t.
\begin{equation}
    x = \bar{\mathcal{X}}(u,w).
\end{equation}    
\end{subequations}

Given this background, the robust optimal control problem with constraints $g(\cdot,\cdot,\cdot):$  $\mathcal{X} \times \mathcal{U} \times \mathcal{W} \rightarrow \mathbb{R}^{n_g}$ that must be satisfied for all admissible values of $w$, can be written as 
\begin{subequations}
    \begin{equation}
        \min_{u \in \mathcal{U}} \, J^*(u)
    \end{equation}
    s.t. $\forall w \in \mathcal{Y}(u)$, $\forall x = \bar{\mathcal{X}}(u,w)$,
    \begin{equation}
        g(x,u,w) \leq 0.
    \end{equation}
\end{subequations}

This problem is very similar to that considered in \cite{wehbeh2024semi}, with the exception of the dependence of the set $\mathcal{Y}$ on $u$. As in \cite{wehbeh2024semi}, we can convert this problem to a more standard semi-infinite programming form through the transformation
\begin{subequations}
\label{eq:robust_optimal_gsip}
    \begin{equation}
        \min_{u \in \mathcal{U}, \gamma \in \Gamma} \, \gamma
    \end{equation}
    s.t. $\forall w \in \mathcal{Y}(u)$,
    \begin{align}
        \gamma - J(\bar{\mathcal{X}}(u,w),u,w) &\leq 0 \\
        g(\bar{\mathcal{X}}(u,w),u,w) &\leq 0
    \end{align}
\end{subequations}
where $\gamma \in \Gamma \subseteq \mathbb{R}$ is a minimal upper bound on the cost that is introduced to enable this reformulation. The resulting problem of \eqref{eq:robust_optimal_gsip} is a generalized SIP, which we reformulate into an existence constrained SIP in Section \ref{sec:gsip_equiv}.

\begin{assum}
    \label{assum:sip_assumptions}
    As in \cite{wehbeh2024semi}, we assume that the set $\mathcal{U}$ is compact, that $\mathcal{W}(x,u)$ is compact $\forall (x,u) \in \mathcal{X} \times  \mathcal{U}$, that the functions $f$, $g$, $h$, and $J$ are continuous in all of their arguments, and that $J(\bar{\mathcal{X}}(u,v),u,w)$ is bounded $\forall u \in \mathcal{U}$, $\forall w \in \mathcal{Y}(u)$, allowing us to also assume that $\Gamma$ is bounded and still obtain an exact solution to \eqref{eq:robust_optimal_gsip}. 
\end{assum}

These assumptions allow us to apply the local reduction method of \cite{wehbeh2024semi} once \eqref{eq:robust_optimal_gsip} is converted into an existence-constrained SIP, as demonstrated in the proof to Theorem~\ref{thrm:convergence} below.


\section{Generalized SIP Equivalence}
\label{sec:gsip_equiv}

We now establish how any generalized SIP can be rewritten as an existence constrained SIP. To do so, we begin by introducing the notation $\bar{u} \coloneqq ( u, \gamma ) \in \bar{\mathcal{U}} \coloneqq \mathcal{U} \times \Gamma$, and define
\begin{equation}
    \bar{g}(\bar{u},w) \coloneqq \begin{bmatrix}
        \gamma - J(\bar{\mathcal{X}}(u,w),u,w) \\
        g(\bar{\mathcal{X}}(u,w),u,w)
    \end{bmatrix}.
\end{equation}
Additionally, we define the objective function $\bar{J}(\cdot):$ $\bar{\mathcal{U}} \rightarrow \mathbb{R}$ such that
\begin{equation}
    \bar{J}(\bar{u}) \coloneqq \gamma
\end{equation}
and introduce the set $\bar{\mathcal{Y}}(\bar{u}) \subseteq \mathbb{R}^{n_w}$ such that
\begin{equation}
\label{eq:uncert_equiv}
    \bar{\mathcal{Y}}(\bar{u}) = \mathcal{Y}(u).
\end{equation}
Using this notation, the problem of \eqref{eq:robust_optimal_gsip} can be rewritten as
\begin{subequations}
\label{eq:gsip_base}
\begin{equation}
    \min_{\bar{u} \in \mathcal{U}} \, \bar{J}(\bar{u})
\end{equation}
s.t. $\forall w \in \bar{\mathcal{Y}}(\bar{u})$,
    \begin{equation}
        \label{eq:gsip_base_constraint}
        \bar{g}(\bar{u},w) \leq 0
    \end{equation}
\end{subequations}
which is the standard form for a generalized SIP. 

Next, we define the set $\bar{\mathcal{W}} \subseteq \mathbb{R}^{n_w}$ such that
\begin{equation}
    \bar{\mathcal{W}} \coloneqq \left\{ w \left| \:
    \exists \bar{u} \in \bar{\mathcal{U}} : \: w \in \bar{\mathcal{Y}}(\bar{u}) 
    \right. \right\}.
\end{equation}
Under this definition, $\bar{\mathcal{W}}$ is equivalent to the union over all of the values of $\bar{u} \in \bar{\mathcal{U}}$ of $\bar{\mathcal{Y}}({\bar{u}})$. Then, we define $\bar{h}(\cdot,\cdot):$ $\bar{\mathcal{U}} \times \mathcal{W} \rightarrow \mathbb{R}^{n_h}$ so that
\begin{equation}
    \label{eq:h_equiv}
    \bar{h}(\bar{u},w) \coloneqq h(\bar{\mathcal{X}}(u,w),u,w)
\end{equation}
and introduce the notation $\bar{g}_i(\bar{u},w)$ to refer to the $i$-th constraint of $\bar{g}(\bar{u},w)$, and $\bar{h}_i(\bar{u},w)$ to the $i$-th constraint of $\bar{h}(\bar{u},w)$. We are now ready to introduce the principal result for this work.

\begin{thrm}
\label{thrm:transformation}
Any generalized SIP of the form of \eqref{eq:gsip_base} can be rewritten as the existence-constrained SIP
\begin{subequations}
\label{eq:existence_sip}
\begin{equation}
    \min_{\bar{u} \in \mathcal{U}} \, \bar{J}(\bar{u})
\end{equation}
s.t. $\forall w \in \bar{\mathcal{W}}$,
    \begin{equation}
    \label{eq:existence_sip_constraint}
    \begin{split}
    \exists \, \lambda \in \Lambda &: \: \ell_1 \max_{i \in \{1,\ldots,n_h\}} \, \left( \bar{h}_i(\bar{u},w) + \epsilon \right) \\[3pt]
    & \qquad + \ell_2 \max_{i \in \{1,\ldots,n_g + 1\}} \, \bar{g}_i(\bar{u},w) \leq 0
    \end{split}
    \end{equation}
\end{subequations}
for some $\epsilon > 0$ and where $\lambda \in \Lambda \subset \mathbb{R}^2$ is a member of the 2-dimensional simplex set with elements $\ell_1$ and $\ell_2$, such that
\begin{subequations}
    \label{eq:lambda_def}
    \begin{align}
        \ell_1 + \ell_2 &= 1 \\
        \ell_1,\, \ell_2 &\geq 0.
    \end{align}
\end{subequations}
The problem of \eqref{eq:existence_sip} will have the same solution $\bar{u}^*$ and associated objective value $\bar{J}(\bar{u}^*)$ as that of \eqref{eq:gsip_base}.
\end{thrm}

\begin{proof}
Since $\bar{\mathcal{W}}$ contains every feasible value of $w$  associated with any choice of $\bar{u}$ that is compatible with the dynamics, we can rewrite the generalized SIP of \eqref{eq:gsip_base} as 
\begin{subequations}
\label{eq:gsip_if}    
\begin{equation}
    \min_{\bar{u} \in \mathcal{U}} \, \bar{J}(\bar{u})
\end{equation}
s.t. $\forall w \in \bar{W}$,
    \begin{equation}
    \label{eq:gsip_if_constraint}
       w \in \bar{\mathcal{Y}}(\bar{u}) \implies \bar{g}(\bar{u},w) \leq 0
    \end{equation}
\end{subequations}
which can be considered a regular SIP because the uncertainty set $\bar{W}$ is independent of the decision variable $\bar{u}$.  

Of course, the constraint of \eqref{eq:gsip_if_constraint} is not trivial to deal with, and the problem of \eqref{eq:gsip_if} is not particularly useful in practice. Therefore, we rewrite the implication of \eqref{eq:gsip_if_constraint} as
\begin{equation}
\label{eq:logical_or_constraint}
    \left[ \neg \left(w \in \bar{\mathcal{Y}}(\bar{u}) \right) \, \right] \lor \left[ \, \bar{g}(\bar{u},w) \leq 0 \vphantom{\bar{\mathcal{Y}}}\, \right].
\end{equation}
Next, we use \eqref{eq:uncertainty_set}, \eqref{eq:uncert_equiv}, and \eqref{eq:h_equiv} to rewrite \eqref{eq:logical_or_constraint} as
\begin{equation}
    \label{eq:function_or_constraint}
    \left[\,\bar{h}(\bar{u},w) < 0 \, \right] \lor \left[ \, \bar{g}(\bar{u},w) \leq 0 \vphantom{\bar{\mathcal{Y}}}\, \right].
\end{equation}
Consequently, since every element $g_i(\bar{u},w)$ of $g(\bar{u},w)$ being negative is equivalent to the largest element being negative, and the same is true for $h(\bar{u},w)$, the constraint of \eqref{eq:function_or_constraint} is equivalent to
\begin{equation}
\label{eq:max_or_constraint}
\begin{split}
        &\left[\,\max_{i \in \{1,\ldots,n_h\}} \, \bar{h}_i(\bar{u},w) < 0 \, \right] \lor \\[3pt]
        &  \qquad\qquad \left[ \, \max_{i \in \{1,\ldots,n_g + 1\}} \, \bar{g}_i(\bar{u},w) \leq 0 \vphantom{\bar{\mathcal{Y}}}\, \right].
\end{split}
\end{equation}
From the definition of an inequality, there must exist an $\epsilon > 0 \in \mathbb{R}$ that allows us to relax the strict inequality into a regular one. Following that, we can convert the logical or ($\lor$) into the equivalent minimum representation as per \cite[Chpt~1]{o2014analysis} and therefore convert \eqref{eq:max_or_constraint} into the inequality
\begin{equation}
\label{eq:min_max_constraint}
\begin{split}
        &\min \left[\,\max_{i \in \{1,\ldots,n_h\}} \, \left( \bar{h}_i(\bar{u},w) + \epsilon \right) \right. , \\[3pt]
        &  \qquad\qquad \left. \max_{i \in \{1,\ldots,n_g + 1\}} \, \bar{g}_i(\bar{u},w)  \, \right] \leq 0. 
\end{split}
\end{equation}

By applying the smoothing approach of \cite{kirjner1998conversion} to the $\min$ operator in \eqref{eq:min_max_constraint}, we recover \eqref{eq:existence_sip_constraint}. 
This implies that the constraint of \eqref{eq:gsip_base_constraint} holds $\forall w \in \bar{\mathcal{Y}}(\bar{u})$ if and only if the constraint of \eqref{eq:existence_sip_constraint} holds $\forall w \in \bar{\mathcal{W}}$. Thus, solutions to \eqref{eq:gsip_base} and \eqref{eq:existence_sip} must satisfy the same constraint. Given that the two problems also share the same objective and decision variable space, their solutions must be identical. The equivalence in their optimal objective values follows trivially.
\end{proof}
\section{Local Reduction on Generalized SIPs}
\label{sec:local_red_gsip}

\subsection{Proof of Convergence}

\begin{thrm}
    \label{thrm:convergence}
    The local reduction algorithm of \cite{wehbeh2024semi} will converge to a (global optimal) solution of \eqref{eq:existence_sip} under Assumptions \ref{assum:f_uniqueness} and \ref{assum:sip_assumptions}. 
\end{thrm}

\begin{proof}
We begin by rewriting $\mathcal{Y}(u)$ as
\begin{equation}
\begin{split}
    \mathcal{Y}(u) &= \left\{ w \left| \: \exists \, x \in \mathcal{X} : \: \begin{array}{r}
        h(x,u,w) \geq 0 \\
        f(x,u,w) = 0
    \end{array}
    \right. \right\} \\
    &= \bigcup_{x \in \mathcal{X}} \mathcal{W}(x,u) \cap \{w\left|\:f(x,u,w) = 0\right.\}.
\end{split}
\end{equation}
We know that the sets $\mathcal{W}(x,u)$ and $\mathcal{U}$ are compact from Assumption~\ref{assum:sip_assumptions}, and that $\{w \,| \:f(x,u,w) = 0\}$ is compact from Assumption~\ref{assum:f_uniqueness}. The set intersection is also compact since it is the intersection of two compact sets. Consequently, from~\cite[Thm~2.5]{michael1951topologies}, $\mathcal{Y}(u)$ must be compact as the union of compact sets over a compact set.

Next, we use \eqref{eq:uncertainty_set} and \eqref{eq:uncert_equiv} to rewrite $\bar{\mathcal{W}}$ as
\begin{equation}
\begin{split}
    \bar{\mathcal{W}} &= \left\{ w \left| \: 
    \exists {u} \in {\mathcal{U}} : \:
    h(\bar{\mathcal{X}}(u,w),u,w) \geq 0
    \right. \right\}  \\
    &= \bigcup_{u \in \mathcal{U}} \left\{ w \left| \:h(\bar{\mathcal{X}}(u,w),u,w) \geq 0 \right. \right\} \\
    &= \bigcup_{u \in \mathcal{U}} \mathcal{Y}(u).
\end{split}
\end{equation}
Applying \cite[Thm~2.5]{michael1951topologies} once again, the uncertainty set $\bar{\mathcal{W}}$ must also be compact since it is likewise a compact union of compact sets. 

From its definition in \eqref{eq:lambda_def}, the set $\Lambda$ is both closed and bounded, and therefore compact by the Heine-Borel theorem. The decision set $\bar{\mathcal{U}}$ is also compact because of the compactness of $\mathcal{U}$ and $\Gamma$ under Assumption \ref{assum:sip_assumptions}. 

We also know that $\bar{J}(\cdot)$ is trivially continuous by definition, and $\bar{g}(\cdot,\cdot,\cdot)$ is continuous as a consequence of the continuity of $g(\cdot,\cdot,\cdot)$ and $J(\cdot,\cdot,\cdot)$ under Assumption \ref{assum:sip_assumptions}. We also know that $\bar{h}(\cdot,\cdot,\cdot)$ is continuous from Assumption \ref{assum:sip_assumptions} and the equivalence of \eqref{eq:h_equiv}. Given that the finite $\max$ operator preserves continuity, we therefore know that the constraint of \eqref{eq:existence_sip_constraint} must be continuous. 

As a result, since all of the functions involved in the definition of \eqref{eq:existence_sip} are continuous, and each of its decision set, uncertainty set, and existence constraint set are compact, the problem of $\eqref{eq:existence_sip}$ can be solved as a series of nested SIPs following the procedure of \cite{wehbeh2024semi}, and will converge to optimality under local reduction as per the results in  \cite{blankenship1976infinitely}. 
\end{proof}

\subsection{Local Reduction Implementation}
\label{sec:local_red_gsip_impl}

Following the algorithm of \cite{wehbeh2024semi}, we solve the existence-constrained problem of \eqref{eq:existence_sip} by converting into a series of minimization and maximization problems. The local reduction approach replaces the continuous constraint set $\bar{\mathcal{W}}$ by the finite scenario approximation $\mathbb{W}_M$ with $M$ different scenarios denoted $w_i$, $i \in \{0,\ldots,M\}$.

These scenarios are found one at a time by maximizing the constraint violation given the current best guess $\bar{u}_k$ for the decision variables. Each constraint $w_{k+1}$ is found as per~\cite{wehbeh2024semi} by solving the optimization problem
\begin{subequations}
    \label{eqn:lr_max}
    \begin{equation}
        w_{k+1} = \argmax_{w \in \bar{\mathcal{Y}}(\bar{u}_k)} \: \sigma
    \end{equation}
s.t. $\forall \lambda \in \Lambda$,
    \begin{equation}
    \label{eqn:lr_max_constraint}
    \begin{split}
        &\sigma - \ell_1 \max_{i \in \{1,\ldots,n_h\}} \, \left( \bar{h}_i(\bar{u}_k,w) + \epsilon \right) \\[3pt]
        & \qquad - \ell_2 \max_{i \in \{1,\ldots,n_g + 1\}} \, \bar{g}_i(\bar{u}_k,w) \leq 0
\end{split}
\end{equation}
\end{subequations}
which is itself a SIP that can be solved using a more standard local reduction approach such as in \cite{blankenship1976infinitely}. 

Once $w_{k+1}$ has been computed, it is added to the previous constraint set $\mathbb{W}_k$ to produce the updated scenario set $\mathbb{W}_{k+1}$. In order to then solve the minimization problem, we begin by defining the new decision variable $L_k = \{ \lambda_1, \ldots,\lambda_k \}$, $\lambda_i \in \Lambda$, $i \in \{1,\ldots,k\}$ to be a collection of $k$ different values of $\lambda$. Then, we obtain $u_{k+1}$ by solving the minimization
\begin{subequations}
\label{eqn:lr_min}
    \begin{equation}
        u_{k+1} = \argmin_{\bar{u} \in \bar{\mathcal{U}}} \min_{L_{k+1} \in \Lambda^{k+1}} \: \bar{J}(\bar{u})
    \end{equation}
    s.t. $\forall \left\{w_i \in \mathbb{W}_{k+1}, \,\lambda_i \in L_{k+1} \right\}$,
    \begin{equation}
    \label{eqn:lr_min_constraint}
    \begin{split}
        &\ell_1 \max_{j \in \{1,\ldots,n_h\}} \, \left( \bar{h}_j(\bar{u},w_i) + \epsilon \right) \\[3pt]
        & \qquad + \ell_2 \max_{j \in \{1,\ldots,n_g + 1\}} \, \bar{g}_j(\bar{u},w_i) \leq 0
    \end{split}
    \end{equation}
\end{subequations}
which is a finitely constrained optimization problem with only $k+1$ constraints. Note, however, that the constraints of \eqref{eqn:lr_max_constraint} and \eqref{eqn:lr_min_constraint} are discontinuous, forcing us to use non-smooth optimization methods to solve \eqref{eqn:lr_max} and \eqref{eqn:lr_min} or to resort to other approaches to smooth out the finite min. 


\section{Simulation Results}
\label{sec:results}

\begin{figure}[b]
    \centering
    \includegraphics[width=0.55\columnwidth]{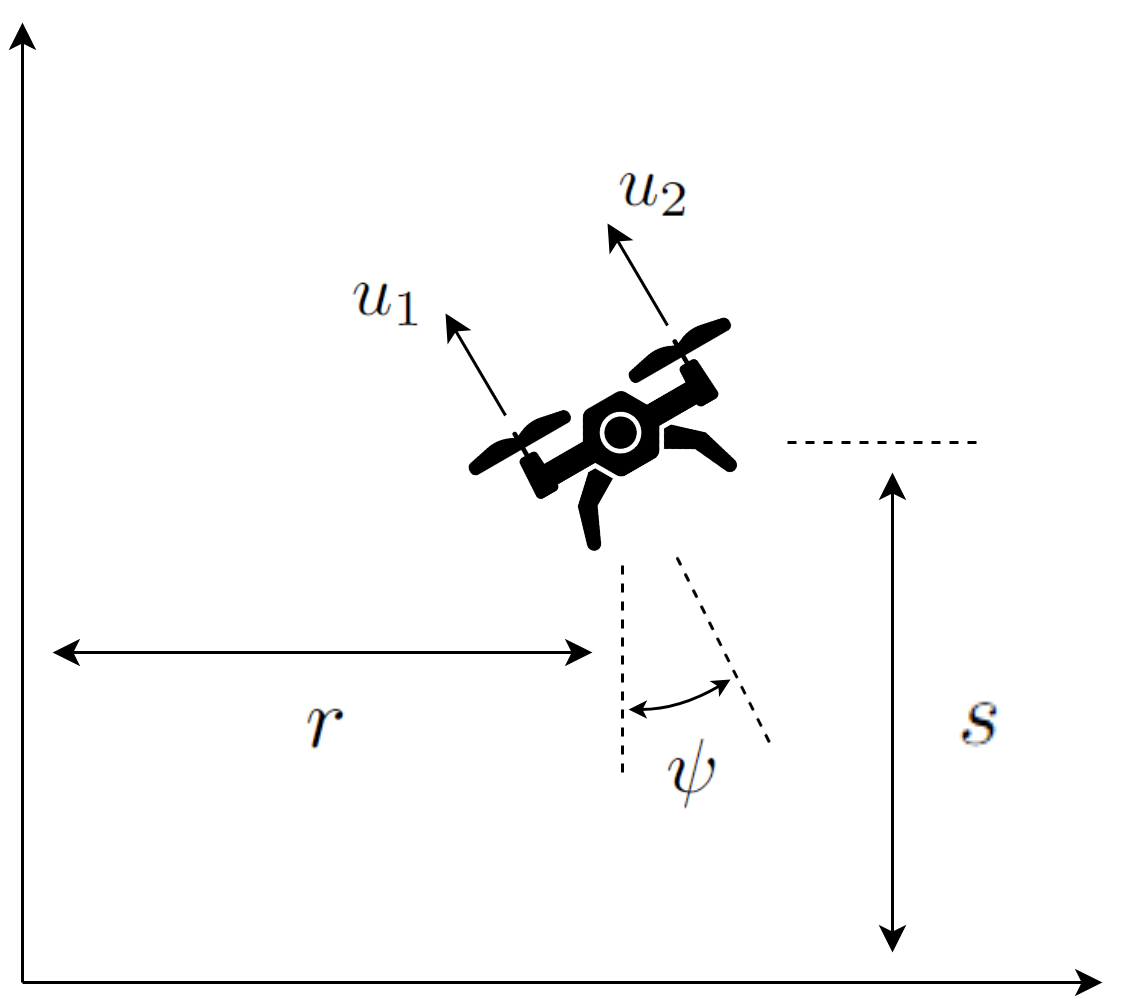}
    \caption{Illustration of the quadrotor's horizontal position ($r$), altitude ($s$), and tilt angle ($\psi$) taken from \cite{wehbeh2024robust}.}
    \label{fig:quad_example}
\end{figure}

We now evaluate the performance of our proposed approach on a simulation of a nonlinear planar quadrotor model with states $[r,\dot{r},s,\dot{s},\psi,\dot{\psi}]$, where $r$ is the quadrotor's horizontal position, $s$ is the quadrotor's height, and $\psi$ is the quadrotor's tilt angle, as illustrated in Figure \ref{fig:quad_example}.

In continuous-time, the quadrotor's positional and angular dynamics are described by the equations
\begin{equation}
    \begin{bmatrix}
        \ddot{r}(t) \\ \ddot{s}(t) \\ \ddot{\psi}(t)
    \end{bmatrix}
    =
    \begin{bmatrix}
        \sin(\psi(t)) \left(u_1(t) + u_2(t) \right) /m \\
        \cos(\psi(t)) \left(u_1(t) + u_2(t) \right) /m - b \\
        \ell (u_1(t) -u_2(t))/ I
    \end{bmatrix}
\end{equation}
where $u_1(t)$ and $u_2(t)$ are the real motor thrusts at time $t$, $m = 0.15$ is the vehicle's mass, $I = 0.00125$ is the moment of inertia, $\ell = 0.1$ is the moment arm for each motor, and $b = 9.81$ is the gravity acting on the system. 

We then obtain the discrete-time dynamics at time step $k$ through the application of a midpoint rule, such that
\begin{equation}
\label{eq:ex_state_dyn}
    \begin{bmatrix}
        {x}_{1,k+1} \\
        {x}_{2,k+1} \\
        {x}_{3,k+1} \\
        {x}_{4,k+1} \\
        {x}_{5,k+1} \\
        {x}_{6,k+1} 
    \end{bmatrix}
    =
    \begin{bmatrix}
        {x}_{1,k} \\
        {x}_{2,k} \\
        {x}_{3,k} \\
        {x}_{4,k} \\
        {x}_{5,k} \\
        {x}_{6,k} 
    \end{bmatrix} +
    T_s
    \begin{bmatrix}
    \frac{\left(x_{2,k} + x_{2,k+1} \right)}{2} \\[3pt]
    \frac{\left(\ddot{r}(kT_s) + \ddot{r}((k+1)T_s) \right)}{2} \\[3pt]
    \frac{\left(x_{4,k} + x_{4,k+1} \right)}{2} \\[3pt]
    \frac{\left(\ddot{s}(kT_s) + \ddot{s}((k+1)T_s) \right)}{2} \\[3pt]
    \frac{\left(x_{6,k} + x_{6,k+1} \right)}{2} \\[3pt]
    \frac{(\ddot{\psi}(kT_s) + \ddot{\psi}((k+1)T_s) )}{2}
    \end{bmatrix}
\end{equation}
where $T_s = 0.2$ is the sampling time, $x_{1,k}$ through $x_{6,k}$ are the discretized states at time $k$ corresponding to the continuous-time values $[r,\dot{r},s,\dot{s},\psi,\dot{\psi}]$, and the discretization is such that $k = 0$ corresponds to $t = 0$. We examine the behavior of the quadrotor over a simulation of $N = 10$ time steps.

In this example, the thrust generated by the quadrotor is uncertain and increases as the thrust commanded increases. The uncertainty is also allowed to vary over time, and is constrained to belong to the set
\begin{equation}
\label{eq:ex_uncert}
 \mathcal{W}_1(v) \coloneqq \left\{ \, w \: \left| \: \: \begin{array}{c}
     -0.05 v_{1,i} \leq w_{1,i} \leq 0.05v_{1,i} \\
     -0.05 v_{2,i} \leq w_{2,i} \leq 0.05v_{2,i} \\
     w_{1,i}\,v_{2,i} = w_{2,i}\,v_{1,i} \\
     i = 1,\ldots,N
\end{array} \right. \right\}   
\end{equation}
where $v \in \mathcal{V} \subset \mathbb{R}^{2N}$ such that $-2 \leq v_{1,i} \leq 2$ and $-2 \leq v_{1,i} \leq 2$ are the reference thrusts for $u_{1,i}$ and $u_{2,i}$ respectively, and
\vspace{-5pt}
\begin{subequations}
\label{eq:ex_contrl_dyn}
\begin{align}
u_{1,i} &= v_{1,i} + w_{1,i} \\
u_{2,i} &= v_{2,i} + w_{2,i}.
\end{align}
\end{subequations}

The objective for this problem is to minimize the function
\begin{equation}
\label{eq:ex_cost}
    J(x) = \sum_{i = 1,\ldots,N} \left( x_{1,i} - r_{\text{ref}}  \right)^2 + \left( x_{3,i} - s_{\text{ref}}  \right)^2
\end{equation}
where $r_{\text{ref}} = 1$ and $s_{\text{ref}} = 2$ and the drone starts from the initial condition $x_0 = [0,0,0,0,0,0]$. We also require that the quadrotor satisfy the constraint
\begin{equation}
\label{eq:ex_constr}
    0 \leq x_{3,i} \leq 2.5 \quad \forall i \in \{1,\ldots,N\}.
\end{equation}

For the rest of this section, we will refer to the problem of minimizing the worst case of the cost of \eqref{eq:ex_cost} subject to the uncertainties of \eqref{eq:ex_uncert}, the dynamics of \eqref{eq:ex_state_dyn} and \eqref{eq:ex_contrl_dyn}, and the constraints of \eqref{eq:ex_constr} as the generalized SIP or GSIP. We will also call the version of the problem reformulated according to the procedure of Section \ref{sec:gsip_equiv} the existence-constrained SIP or ESIP. We compare our solution to the relaxed SIP, referred to henceforth as SIP 1, with the uncertainty set
\begin{equation}
    \bar{\mathcal{W}}_1 \coloneqq \bigcup_{v \in \mathcal{V}} \mathcal{W}_1(v) = \left\{ \, w \: \left| \: \: \begin{array}{c}
     -0.1 \leq w_{1,i} \leq 0.1 \\
     -0.1 \leq w_{2,i} \leq 0.1 \\
     - w_{1,i} \, w_{2,i} \leq 0 \\
     i = 1,\ldots,N
\end{array} \right. \right\}
\end{equation}
which does not depend on the decision variables. We also compare our approach to SIP 2 with the uncertainty set
\begin{equation}
    \bar{W}_2 \coloneqq\left\{ \, w \: \left| \: \: \begin{array}{c}
     -0.1 \leq w_{1,i} \leq 0.1 \\
     -0.1 \leq w_{2,i} \leq 0.1 \\
    w_{1,i} = w_{2,i} \\
     i = 1,\ldots,N
\end{array} \right. \right\}
\end{equation}
which is a less conservative approximation of $\mathcal{W}_1(v)$ that does not satisfy $\mathcal{W}_1(v) \subset \bar{\mathcal{W}}_2$ $\forall v \in \mathcal{V}$. Finally, we also consider the exact reformulation of the GSIP (RSIP) using $w' \in \mathcal{W}^{\,\prime}_1$ such that
\begin{equation}
    \mathcal{W}_1^{\,\prime} \coloneqq \left\{ \, w^{\prime} \: \left| \: \: \begin{array}{c}
     -0.05 \leq w_i^{\prime} \leq 0.05 \\[3pt]
     i = 1,\ldots,N
\end{array} \right. \right\}
\end{equation}
\begin{subequations}
\begin{align}
u_{1,i} &= v_{1,i}\,(1 + w_i^{\prime})  \\
u_{2,i} &= v_{2,i}\,(1 + w_i^{\prime}).
\end{align}
\end{subequations}

\begin{figure}[tb]
    \centering
    \includegraphics[width=0.7\columnwidth]{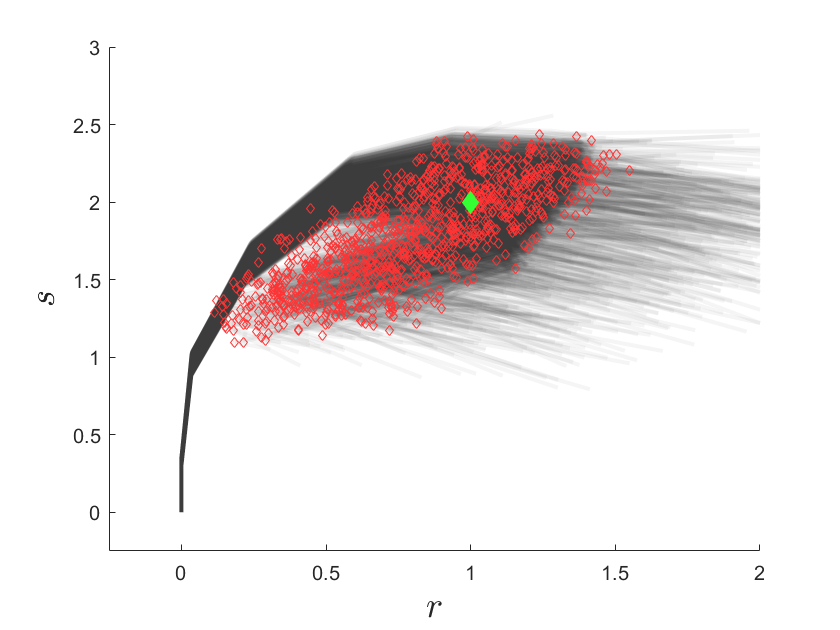}
    \vspace{-5pt}
    \caption{Quadrotor trajectories (in gray) for the solution obtained from the ESIP or RSIP. Final positions are shown in red and the target in green.}
    \label{fig:esip_sol}    
    \vspace{-5pt}
\end{figure}

Note that while the problems considered here all compute open-loop policy, our methods can be generalized to compute optimal feedback gains as was done in \cite{wehbeh2024robust}.

\begin{table}[tb]
\centering
\caption{Performance results for key metrics across 10000 simulation runs on the four problems considered.}
\label{tab:results}
\begin{tabular}{@{}lllll@{}}
Problem               & ESIP  & SIP 1   & SIP 2 & RSIP  \\ \toprule
$\gamma$              & 16.04 & No sol. & 19.05 & 16.04 \\ \midrule
Sol. Scen.            & 12    & No sol. & 7     & 12    \\ \midrule
Avg. Cost             & 12.09 & No sol. & 22.56 & 12.09 \\ \midrule
Worst Cost            & 15.34 & No sol. & 88.79 & 15.34 \\ \midrule
Constr. Viol.         & 0     & No sol. & 455   & 0    
\end{tabular}
\vspace{-15pt}
\end{table}

Each problem is modeled in Julia using the JuMP package~\cite{Lubin2023} according to the procedure described in Section \ref{sec:local_red_gsip_impl} or in \cite{wehbeh2024semi}. The optimization problems obtained are then solved by using the Juniper optimizer \cite{juniper} to deal with the non-smooth constraints and Ipopt \cite{wachter2006implementation} to handle the nonlinear problem. After open-loop control policies were obtained for each problem, 10000 different uncertainty realizations were sampled from \eqref{eq:ex_uncert} for each control trajectory to evaluate their performance. As expected, the ESIP and RSIP return the same solution. SIP 1 was found to be infeasible and could not be solved. The results are summarized in Table~\ref{tab:results}.

\begin{figure}[tb]
    \centering
    \includegraphics[width=0.7\columnwidth]{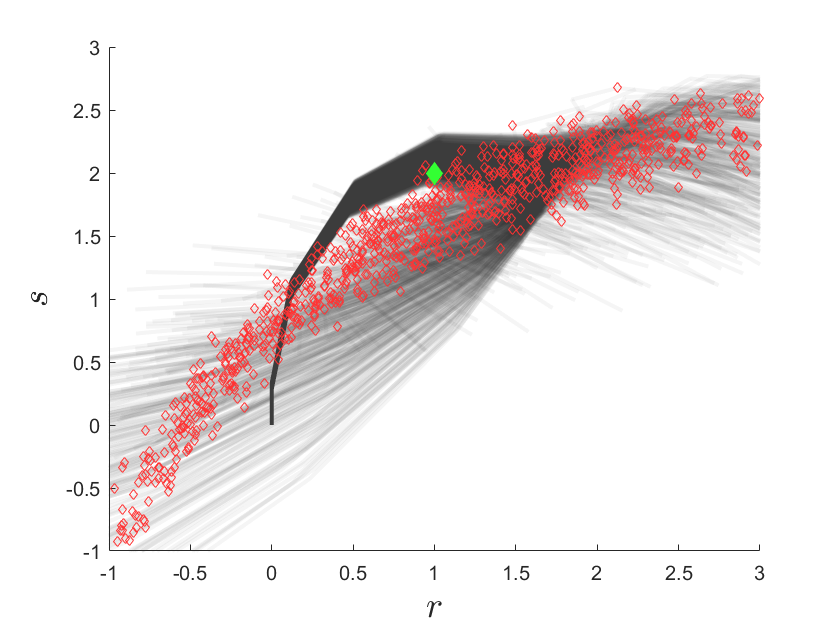}
    \vspace{-5pt}
    \caption{Quadrotor trajectories for the solution obtained from SIP 2. Colors are as in Figure \ref{fig:sip2_sol}.}
    \vspace{-10pt}
    \label{fig:sip2_sol}
\end{figure}

Since the ESIP and RSIP yield the same solution, their entries in Table \ref{tab:results} are identical. Both approaches converge after 12 local reduction iterations, produce a robust solution that leads to no constraint violations, and compute a cost bound $\gamma$ that is greater than the largest cost seen in simulation. As seen in Figure \ref{fig:esip_sol}, most trajectories converge to the proximity of the target despite the uncertainty. The solution of SIP 2, meanwhile, fails to properly model the uncertainty, and performs worse in practice than the solution predicts. A large number of trajectories lead to constraint violations or stray far from the equilibrium as seen in Figure \ref{fig:sip2_sol}.

\section{Conclusions}
\label{sec:conclusion}
In this work, we introduced a novel method to transform any generalized SIP into an existence-constrained SIP, and explained how such problems can be solved in practice. We then demonstrated this on a planar quadrotor control example, where our method yielded the same results as an exact reformulation of the generalized SIP and outperformed more conservative transformations. This implies that our method can be used on problems where exact reformulations are much harder to define and still produce exact solutions. 



\section*{ACKNOWLEDGMENTS}

This work was funded by the Natural Sciences and Engineering Research Council of Canada through a PGS D grant. 


\bibliographystyle{IEEEtran}
\bibliography{IEEEabrv,references.bib}

\end{document}